\documentclass[a4paper,11pt]{article}

\usepackage{amsmath,amsfonts,amscd,amssymb,amsthm}
\usepackage[dvips]{graphics}

\newtheorem{theorem}{Theorem}[section]
\newtheorem*{thm}{Theorem}

\newtheorem{lemme}[theorem]{Lemme}

\newtheorem{definition}{Definition}[section]
\newtheorem{coro}{Corollary}[section]

\newtheorem{prop}{Proposition}[section]

\newcommand{\dd}{\text{d}}
\newcommand{\OO}{\mathcal{O}}
\newcommand{\DD}{\mathcal{D}}
\newcommand{\XX}{\mathcal{X}}
\newcommand{\FF}{\mathcal{F}}

\newcommand{\MM}{\mathcal{M}}
\newcommand{\UU}{\mathcal{U}}

\newcommand{\coup}[1]{\left( #1 \right)}

\newcommand{\carre}{
\begin{flushright}
$\square$
\end{flushright}
\medskip}

\newcommand{\fraction}[2]{\displaystyle\frac{#1}{#2}}
\newenvironment{demo}{\noindent {\it Proof:}}{\carre}

\begin{document} 

\title{Rigidity for dicritical germ of foliation in $\mathbb{C}^2$.}

\maketitle

\abstract{We study some kind of rigidity property for dicritical foliation in $\mathbb{C}^2$ thanks to a new meromorphic function associated to any dicritical foliation playing the role of equation of separatrix . To be more specific, we compute infinitesimal obstructions for realizing the holonomy pseudo-group of a dicritical foliation $\FF$ on any source manifold of a blowing-up morphism with same dual tree as the desingularization of $\FF$. Using two examples, we remark that these obstructions come from some analytical invariants shared by \emph{the space of leaves} and the ambient space. This situation never occurs in the non-dicritical case.}

\section*{Introduction and main statements}
Considering the problem of moduli for a germ of singular holomorphic foliation $\FF$ in $\mathbb{C}^2$ leads naturally to point out some topological and analytical invariants. The invariants of first kind come from the reduction of singularities $E:(\MM,\DD)\rightarrow (\mathbb{C}^2,0)$ of the foliation: a combinatorial invariant given by the topological class of the manifold $\MM$ and a analytical one given by the analytical class of the manifold $\MM$, called the resolution space of $\widehat{\FF}$. The invariants of second kind are more related to the foliation itself: the collection of projective holonomy representations defined over each component of the divisor $\DD$ and Dulac applications, in other words \emph{the holonomy pseudo-group}. A natural problem is to know whether \emph{coherent} data of above invariants correspond to a concrete foliation. 

\noindent A germ of formal foliation $\widehat{\FF}$ in $\mathbb{C}^2$ is the data of a formal germ of $1$-form at $0\in\mathbb{C}^2$
\begin{equation}\label{formalform}
\omega={a}(x,y)\dd x + {b}(x,y)\dd y
\end{equation}
up to a unity and a conjugacy with $a,b\in\mathbb{C}[[x,y]]$. When $\omega$ is convergent, $\widehat{\FF}$ is said to be convergent and will be denoted simply by $\FF$. A \emph{formal separatrix} of $\widehat{\FF}$ is a formal irreducible curve $\{f=0\}$ given by a formal germ of reduced function $f$ such that $f$ divides $\dd f \wedge  \omega$ in $\mathbb{C}[[x,y]]$. If $\FF$ is convergent, a separatrix is nothing but a germ of analytical curve which is a formal separatrix of $\FF$ seen as a formal foliation. In \cite{separatrice}, C. Camacho and P. Sad proved that any singular germ of foliation in $(\mathbb{C}^2,0)$ admits at least one separatrix. A foliation is said \emph{dicritical} when it has not a infinite number of separatrix. 

\noindent Let us denote by $E:(\MM,\DD)\rightarrow (\mathbb{C}^2,0)$ the reduction process (\cite{reduction}\cite{MM}) of $\widehat{\FF}$ where $\DD$ refers to the exceptional divisor $E^{-1}(0)$. The singularities of $E^*\widehat{\FF}$ are \emph{reduced}, which means that, in some local coordinates, they belong to the following list
\begin{enumerate}
\item  $\lambda x\dd y + y\dd x + \cdots $ terms of higher order,\quad $\lambda \not\in \mathbb{Q}^-$,
\item $x\dd y +\ldots$ terms of higher order.
\end{enumerate}
In the second case, the singularity is called \emph{saddle-node}. The formal normal forms for such singularity are given in \cite{MatSal}: there exist some coordinates so that the singularity is given by the following $1$-form
$$ (\zeta x^p-p)y\dd x+x^{p+1}\dd y,\quad p\in \mathbb{N}^*,\zeta\in \mathbb{C}.$$
The invariant curve $\{x=0\}$ is called the \emph{strong invariant curve} and $\{y=0\}$ the \emph{weak one}. If the germ of the divisor $\DD$ is the weak invariant curve, then the singularity is called \emph{tangent saddle-node}. We recall that $\FF$ is said to be of \emph{second kind} when $\widehat{\FF}$ is non-dicritical and when none singularity of $E^*\FF$ is a tangent saddle-node. One have same definition in the convergent context with obvious transposition.

\medskip
\noindent In this article, we study the extended notion of \emph{foliation of second kind}. The definition is the same as \emph{foliation of second kind} defined in \cite{MatQuasi} but, of course, we do not require the foliation to be non-dicritical. The construction of \emph{a balanced equation} of the separatrix allows us to study this class of foliations. In the non-dicritical case, the balanced equation coincide with an equation of the finite number of separatrix and we prove a criterion for a foliation to be of second kind, which is a standard criterion in the non-dicritical case \cite{MatQuasi}. 

\noindent Let us define $\widehat{\XX}_{\widehat{F}=0}$ the sheaf of base $\DD$ whose fibre is the space of germs of vector field tangent to the total transform of the zero of a balanced equation $\widehat{F}$. Let $\widehat{\XX}_{\widehat{\FF}}$ be the sub-sheaf of vector field tangent to the foliation $\widehat{\FF}$. The sheaf $\widehat{\mathfrak{I}}_{ \widehat{ F}=\infty}$ refers to the germs of functions vanishing along the pole of $\widehat{F}$. The first 
\begin{thm}\label{CritereEsp}
Let $\widehat{F}$ be a balanced equation. The following properties are equivalent:
\begin{enumerate}
\item $\widehat{\FF}$ is of second kind.
\item $\nu_{0}(\widehat{F})=\nu_0(\widehat{\FF})+1$.
\item The sequence of sheaves
$$0\longrightarrow \widehat{\XX}_{\widehat{\FF}} \longrightarrow \widehat{\XX}_{\widehat{F}=0} \xrightarrow{E^*\frac{\widehat{\omega}}{\widehat{F}}(\cdot)}\widehat{\mathfrak{I}}_{ \widehat{F}=\infty}\longrightarrow 0$$ is exact.
\end{enumerate}
\end{thm}
\noindent Thanks to the above result, we compute the infinitesimal obstructions for the following \emph{problem of existence}:
\begin{quote}
can any deformation of the resolution space of $\widehat{\FF}$ be followed by an isoholonomical deformation of $\widehat{\FF}$ ?
\end{quote}
 This question is the local version of an harder problem: can any manifold topologically equivalent to the resolution space of $\widehat{\FF}$ carry a foliation with same holonomy groupo\"id as $\widehat{\FF}$. In the second part of this article, we will give a precise statement for this kind of problem. Finally, we are able to prove the following result:
\begin{thm}For a generic dicritical formal foliation $\widehat{\FF}$, 
there exists an analytical topologically trivial deformation of the space of resolution of $\widehat{\FF}$, which cannot carry any isoholonomical deformation of $\widehat{\FF}$.
\end{thm}

\section{Multiplicity of a formal dicritical foliation.}

Our first aim is to establish a formula, which links the multiplicity of $\nu_0(\widehat{\FF})$ and some invariants produced by the desingularization of $\widehat{\FF}$. 

\noindent\emph{The valence $v(D)$} of an irreducible component $D$ of $\DD$ is the number of irreducible components of $\DD$, which intersect $D$ . In the same way, the integer $v_{\bar{d}}(D)$ refers to \emph{the non-dicritical valence}, which means the number of non-dicritical components, which intersect $D$. Let us denote by $\mathfrak{I}_D$ the sheaf over $\DD$ of germs of holomorphic functions, which vanish along the component $D$ of $\DD$. Let us denote by $\mathfrak{M}$ the sheaf over $\DD$ of $\OO$-module generated by the functions $E^*h$ with $h\in\mathcal{O}_2$ and $h(0)=0$. By induction on the number of blowing-up in $E$, one can prove that there exist numbers $\nu(D)$, called \emph{the multiplicity of $D$} such that, one has the following decomposition
$$\mathfrak{M}=\prod_{D\in \text{Comp}(\DD)}\mathfrak{I}_D^{\nu(D)}.$$
The integer $\nu(D)$ is also the multiplicity of a curve whose strict transform is smooth and attached to a regular point of $D$.

\medskip
\noindent Let us consider the following definition introduced by C. Hertling in \cite{hertlingfor}.
\begin{definition} Let $\widehat{\FF}$  be a germ of formal foliation given by a $1$-form 
$$\omega={a}(x,y)\dd x + {b}(x,y)\dd y$$
\noindent \begin{enumerate}
\item Let $(\widehat{S},p)$ a germ of smooth formal invariant curve. If, in some coordinates, $\widehat{S}$ is the curve $\{y=0\}$ and $p$ the point $(0,0)$, then the integer $\textup{ord}_0b(x,0)$ is called the indice of $\widehat{\FF}$ at $p$ with the respect of $\widehat{S}$ and is denoted by $\textup{Ind}(\widehat{\FF},\widehat{S},0)$.

\item Let $(\widehat{S},p)$ a germ of smooth formal non-invariant curve. If, in some coordinates, $\widehat{S}$ is the curve $\{y=0\}$ and $p$ the point $(0,0)$, then the integer $\textup{ord}_0a(x,0)$ is called the tangency order of $\widehat{\FF}$ with the respect of $\widehat{S}$ and is denoted $\textup{Tan}(\widehat{\FF},\widehat{S},p)$. 
\end{enumerate} 
\end{definition}
\noindent  In \cite{hertlingfor}, one can find the following formula, which is an extension of \cite{camacho} when there are dicritical components:
\begin{prop}[\cite{hertlingfor}]\label{Hertlingformule}The multiplicity of $\widehat{\FF}$ satisfies the equality
$$\nu_0(\widehat{\FF})+1=\sum_{D\in \textup{Comp}(d)} \nu(D)\rho(D)$$
where
\begin{enumerate}
\item if $D$ is non-dicritical,
$\rho(D)=-v_{\bar{d}}(D)+\sum_{q\in D} \textup{Ind}(E^*\widehat{\FF},D,q)$.
\item if $D$ is dicritical,
$\rho(D)=2-v_{\bar{d}}(D)+\sum_{q\in D} \textup{Tan}(E^*\widehat{\FF},D,q)$.
\end{enumerate}
\end{prop}
\noindent One has to notice that this formula is true even if the application $E$ is not the exactly the reduction of singularities but any morphisms build over $E$ by a succession of elementary blowing-up.  

\noindent We are going to define a notion of \emph{balanced equation} of separatrix of $\widehat{\FF}$. A separatrix of $\widehat{\FF}$ is said to be \emph{isolated} if its stricted transform with the respect of $E$ is attached to a non-dicritical component. When $D$ is dicritical, one call \emph {pencil of $D$} the set of invariant curves, which strict transform is attached to $D$. 
\begin{definition}
A complete system of separatrix is the union of two germs of curve $Z\cup P$ where
\begin{enumerate}
\item $Z$ is the union of isolated separatrix and, for any dicritical component $D$ with valence smaller than $2$, $2-v(D)$ curves of the pencil of $D$.
\item $P$ is the union of $v(D)-2$ curves of the pencil of any dicritical component $D$ with valence bigger then $3$.
\end{enumerate} 

\noindent A balanced equation of separatrix is a germ of meromorphic function  whose zeros and poles are the blown down at the origin of $\mathbb{C}^2$ of respectively the sets $Z$ and $P$. 
\end{definition}

\medskip
\noindent Let us denote by $\mathcal{SNT}(\widehat{\FF})$ the set of singularity of $E^*\widehat{\FF}$, which are tangent saddle-node. Let $\widehat{F}$ be any balanced equation.
\begin{prop}\label{relationmultiplicite} We have the following equality
$$\nu_0(\widehat{F})=\nu_0(\widehat{\FF})+1+\sum_{s\in \mathcal{SNT}(\widehat{\FF})} \sum_{D\in V(s)}\nu(D)\coup{\textup{Ind}(E^*\widehat{\FF},D,s)-1}$$
where $V(s)$ refers to the set of irreducible components containing the point $s$. In particular, the multiplicity of $\widehat{F}$ does not depend on the choice of $\widehat{F}$.
\end{prop}
\begin{demo} Let us write $\widehat{F}=\fraction{\widehat{N}}{\widehat{P}}$ where $\widehat{N}$ and $\widehat{P}$ are holomorphic. The foliations $d \widehat{N}$ and $d \widehat{P}$ are desingularized by the morphism $E$. By applying (\ref{Hertlingformule}) to $\widehat{\FF}$, $d \widehat{N}$ and $d \widehat{P}$, we get the next relations: 
\begin{eqnarray*}
\nu_0(\widehat{\FF})+1&=&\displaystyle\sum_{D\in \textup{Comp}(d)} \nu(D)\rho(D), \\
\nu_0(d \widehat{N})+1&=&\displaystyle\sum_{D\in \textup{Comp}(d)} \nu_{\widehat{N}}(D)\rho_{\widehat{N}}(D), \\
\nu_0(d \widehat{P})+1&=&\displaystyle\sum_{D\in \textup{Comp}(d)} \nu_{\widehat{P}}(D)\rho_{\widehat{P}}(D). \\
\end{eqnarray*}
As components multiplicity depend only on $E$, we have for any component $D$
$$\nu(D)=\nu_{\widehat{N}}(D)=\nu_{\widehat{P}}(D).$$
Now, if $D$ is non-dicritical for $\FF$, we find
$$\rho(D)=\text{Iso}(D)+\sum_{q\in D\cap\mathcal{SNT}(\widehat{\FF})} (\textup{Ind}(E^*\widehat{\FF},D,q)-1),$$
where $\text{Iso}(D)$ is the number of isolated separatrix attached to $D$.
If $D$ is dicritical, $\rho(D)$ is equal to $2-v(D)$. 
\noindent By definition, if $D$ is non-dicritical, the numbers $\rho_{\widehat{N}}(D)$ and $\text{Iso}(D)$ are equal. Moreover, if $D$ is dicritical with valence smaller than $2$, $\rho_{\widehat{N}}(D)$ is equal to $2-v(D)=\rho(D)$. On any other component, the foliation $\{\dd \widehat{N}=0\}$ doesn't have any isolated separatrix and the integer $\rho_{\widehat{N}}(D)$ vanishes. Hence, we have the relation
\begin{multline}\label{susueq}
\nu_0(\widehat{N})=\nu_0(\dd \widehat{N})+1=\\ \sum_{{\tiny \begin{array}{l}D\in \textup{Comp}(d) \\ D\text{ non-dicritical or} \\ D\text{ dicritical and } v(D)\leq 2\end{array}}} \nu(D)\rho(D)\quad +\hspace{1.5cm}\\ \sum_{{\tiny \begin{array}{l}D\in \textup{Comp}(d) \\ D\text{ non-dicritique} \end{array}}} \nu(D)\sum_{q\in D\cap\mathcal{SNT}(\widehat{\FF})} (\textup{Ind}(E^*\widehat{\FF},D,q)-1).
\end{multline}
As $E^*\widehat{\FF}$ is reduced, any point $q$ in a dicritical component verifies the relation
  $$\text{Tan}(E^* \widehat{\FF},D,q)=0.$$ Hence, for any dicritical components with valence greater than $3$, we have $\rho_{\widehat{P}}(D)=-\rho(D)$.
In any other case, the integer $\rho_{\widehat{P}}(D)$ is zero. Therefore, we have
\begin{equation}\label{susueq2}
\nu_0(\widehat{P})=\nu_0(d \widehat{P})+1=-\sum_{{\tiny \begin{array}{l}D\in \textup{Comp}(d)\\ D\text{ dicritical and } v(D)\geq 3\end{array}}} \nu(D)\rho(D).
\end{equation}
The proposition is the combination of relations (\ref{susueq}) and (\ref{susueq2}).
\end{demo}

\subsection{Dicritical foliation of second kind}

\noindent The notion of balanced equation comes naturally with a notion of dicritical foliation of second kind. Known properties of second kind class will have analogues in our context.
\begin{definition}
$\widehat{\FF}$ is said to be of second kind when none singularity of $E^{*}\widehat{\FF}$ is a tangent saddle-node. When $\FF$ is convergent, $\FF$ is said to be of second kind when it is as formal foliation.
\end{definition}
\noindent Of course, we do not require here $\widehat{\FF}$ to be non-dicritical. The following corollary of (\ref{relationmultiplicite}), which is  shows how these foliations naturally generalized the class of non-dicritical foliation of second kind:
\begin{prop}\label{critere2esp}
Let $\widehat{F}$ be a balanced equation. Then $\widehat{\FF}$ is of second kind if and only if $\nu_{0}(\widehat{F})=\nu_0(\widehat{\FF})+1$
\end{prop}
\noindent This property is a clear analogue of a property for non-dicritical foliation of second kind established in \cite{MatSal}. 

\section{Analysis of dicritical obstructions.} 

 In \cite{moi}, we prove the following result, which can be summary by saying that analytical invariants of the holonomy pseudo-group and analytical invariants of the ambient space are disconnected.
\begin{theorem}[\cite{moi}]\label{superthe}
Let $\FF_0$ be a singular formal {\bf non-dicritical} foliation of second kind at $0\in\mathbb{C}^2$ and $E_0:\MM\rightarrow \mathbb{C}^2$ its desingularization. For any blowing-up process $E_1$ topologically equivalent to $E_0$, there exists a foliation $\FF_1$ at $0\in \mathbb{C}^2$ linked to $\FF_0$ by an isoholonomical deformation such that the desingularization of $\FF_1$ is exactly $E_1$.
\end{theorem}
\noindent For a precise definition of \emph{isoholonomical deformation}, we refer to \cite{univ}; but basically, these are topologically trivial deformations along which the holonomy pseudo-group is constant. The main tool for the proof of the above theorem is the equivalent result at the infinitesimal level, which is a trivial consequence of the existence of the following exact sequence of sheaves
\begin{equation}\label{superseq}
 0\rightarrow \XX_\FF \rightarrow \XX \rightarrow \OO_\MM\rightarrow 0.
\end{equation}
Here, all the sheaves have $\DD$ for base. The fibre of $\XX$ is the space of vector fields tangent to the total transform by $E$ of the separatrix of $\FF$ at the origin of $\mathbb{C}^2$. The sheaf $\XX_\FF$ is the sub-sheaf of $\XX$ whose fibre is the space of vector field tangent to the foliation $E^*\FF$. These two sheaves correspond respectively to the modular space $H^1(\XX)$ of infinitesimal deformation of the ambient space and to the modular space $H^1(\XX_\FF)$ of infinitesimal isoholonomical deformation \cite{univ}. Finally, $\OO_\MM$ refers to the restriction at $\DD$ of the structural sheaf of $\MM$. Since the first cohomology group of $\OO_MM$ is trivial, one has the following sequence
$$H^1(\XX_\FF)\rightarrow H^1(\XX)\rightarrow 0$$
which is the starting point of the proof of (\ref{superthe}).
\subsection{Infinitesimal obstructions.}\label{contreexmero}

In this section, we establish a equivalent of the sequence (\ref{superseq}) in order to compute the formal infinitesimal obstructions, which prevent the theorem (\ref{superthe}) from being true for dicritical foliations.

\begin{prop}\label{CritereEspExa}
The following propositions are equivalent:
\begin{enumerate}
\item $\FF$ is second kind.
\item The sequence of sheaves
$$0\longrightarrow \widehat{\XX}_{\widehat{\FF}} \longrightarrow \widehat{\XX}_{\widehat{F}=0} \xrightarrow{E^*\frac{\widehat{\omega}}{\widehat{F}}(\cdot)}\widehat{\mathfrak{I}}_{ \widehat{F}=\infty}\longrightarrow 0$$
is exact.
\end{enumerate}
\end{prop}
\noindent First, we show that the order of multiplicity of the blown-up balanced equation along any irreducible component of the divisor behave well with respect to the multiplicity of the foliation.

\begin{lemme}\label{muldic}
For any component $D$, we have the following alternative:
\begin{enumerate}
\item if $D$ is non-dicritical, $\nu_D(\widehat{F})=\nu_D(\widehat{\FF})+1$,
\item if $D$ is dicritical, $\nu_D(\widehat{F})=\nu_D(\widehat{\FF})$,
\end{enumerate}
where $\nu_{D}(*)$ refers to the multiplicity of the blown-up in a generic point of $D$.
\end{lemme}
\begin{demo} The proof is an induction on the height of the component $D$ in the blowing-up process. At height $1$, one can see that: the integers $\nu_{D_0}(F)$ and $\nu_0(F)$ are equal; if $D_0$ is non-dicritical, $\nu_{D_0}(\widehat{\FF})=\nu_0(\widehat{\FF})$ else $\nu_{D_0}(\widehat{\FF})=\nu_0(\widehat{\FF})+1$. Hence, the lemma is the proposition (\ref{CritereEsp}). For the induction, we look at the blowing-up process at height $i$ and consider a component $D$ of $\DD^{i+1}$ obtained by blowing-up of $c\in \DD^i$. Let $\widehat{F}^i$ by the divided blown-up of $\widehat{F}$ by $E^{*i}$. We have the following relations:
\begin{eqnarray}
\nu_D(\widehat{\FF})&=&\nu_c(E^{i*}\widehat{\FF})+\displaystyle\sum_{D_c\in\text{V}(c)}\nu_{D_c}(\widehat{\FF}) +\epsilon(D)\label{eqq1}, \\
\nu_D(\widehat{F})&=&\nu_c(\widehat{F}^i)+\displaystyle\sum_{D_c\in\text{V}(c)}\nu_{D_c}(\widehat{F}),\label{eqq2}
\end{eqnarray}
where $\epsilon(D)$ is $0$ when $D$ is non-dicritical, $1$ else. From now on, one has to look at each different cases. We recall that $V(s)$ refers to the set of irreducible components of $\DD$ that contain $s$.
\begin{enumerate}
\item {\bf $V(c)$ consists of one component $D_0$.}
\begin{enumerate}
\item {\bf $D_0$ is non-dicritical}: let $\widehat{F}_c$ the germ of meromorphic function near the point $c$ product of $\widehat{F}^i$ and of a germ of equation of $D_0$. By definition, $\widehat{F}_c$ is a balanced equation for $E^{i*}\widehat{\FF}$. In view of (\ref{CritereEsp}), we have
\begin{equation}\label{eq1}
\nu_c(\widehat{F}_c)=\nu_c(E^{i*}\widehat{\FF})+1.
\end{equation} 
Now, the above construction ensures the equality
\begin{equation}\label{eq2}
\nu_c(\widehat{F}_c)=1+\nu_c(\widehat{F}^i_c).
\end{equation}
Moreover, the induction hypothesis shows the relation
\begin{equation}\label{eq3}
\nu_{D_0}(\widehat{F})=\nu_{D_0}(\widehat{\FF})+1.
\end{equation} 
The association of (\ref{eqq1}), (\ref{eqq2}), (\ref{eq1}), (\ref{eq2}) and (\ref{eq3}) gives the checked result for $D$.
$$\nu_D(\widehat{F})=\nu_D(\widehat{\FF})+1-\epsilon(D).$$
\item {\bf $D_0$ is dicritical}: in that case, one has to choose for $\widehat{F}_c$ the germ $\widehat{F}^i_c$. Once again, $\widehat{F}_c$ is a balanced equation for  $E^{i*}\widehat{\FF}$. Hence,
$\nu_c(\widehat{F}_c)$ is equal to $\nu_c(E^{i*}\widehat{\FF})+1$.
Therefore, one gets the relation $\nu_c(\widehat{F}_c)=\nu_c(\widehat{F}^i_c)$.
Under the induction hypothesis, $\nu_{D_0}(\widehat{F})$ and $\nu_{D_0}(\widehat{\FF})$  are equal.
The association of previous relation ensures the equality 
$$\nu_D(\widehat{F})=\nu_D(\widehat{\FF})+1-\epsilon(D),$$
$D_0$ being dicritical, $\epsilon(D)$ vanishes. 
\end{enumerate}
\item{\bf $V(c)$ consists of two components $D_0$ and $D_1$.}
\begin{enumerate}
\item {\bf $D_0$ and $D_1$ are non-dicritical}: let $\widehat{F}_c$ the product of $\widehat{F}^i$ and of a germ of equation for $D_0\cup D_1$. Since the function $\widehat{F}_c$ is a balanced equation for $E^{i*}\widehat{\FF}$, 
$\nu_c(\widehat{F}_c)$ is equal to $\nu_c(E^{i*}\widehat{\FF})+1$. Now, in view of our construction, we have the equality $\nu_c(\widehat{F}_c)=2+\nu_c(\widehat{F}^i_c)$. Now, for any $D_c\in V(c)$,
$\nu_{D_c}(\widehat{F})=\nu_{D_c}(\widehat{\FF})+1$. 
These equalities ensures that
$$\nu_D(\widehat{F})=\nu_D(\widehat{\FF})+1-\epsilon(D).$$
\item{\bf $D_0$ is dicritical and $D_1$ is non-dicritical}: this case can be treated in the same way.
\end{enumerate} 
\end{enumerate}
\end{demo}

\begin{demo}{\it\  (\ref{CritereEspExa})\ }
Let us suppose $\FF$ of second kind and let $c$ be any point of $\DD$. The kernel of $E^*\frac{\widehat{\omega}}{\widehat{F}}(\cdot)$ clearly intersects $\widehat{\XX}$ along the sheaf $\widehat{\XX}_{\widehat{\FF}}$. Now, we compute the co-kernel over $c$. In each following case, we give a local expression of $E^{*}\frac{\widehat{\omega}}{\widehat{F}}$ in adapted coordinates. The existence of such coordinates is proved in \cite{MatSal} and uses the second kind hypothesis. We give also a solution $X$ in $\widehat{\XX}_{\widehat{F}=0}$ of the equation $E^*\frac{\widehat{\omega}}{\widehat{F}}(X)=g$ for $g$ in $\widehat{\mathfrak{I}}_{ \widehat{F}=\infty}$:
\begin{enumerate} \item  $c$ is a regular point of a non-dicritical component, which is neither a zero nor a pole of $\widehat{F}$:
$$E^*\frac{\widehat{\omega}}{\widehat{F}}=\fraction{u}{x}\dd x,\quad X=\fraction{g}{u}x\fraction{\partial}{\partial x} ,\quad \coup{\widehat{\mathfrak{I}}_{ \widehat{F}=\infty}}_c=\OO_c.$$
\item $c$ is a regular point of a dicritical component, which is neither a zero nor a pole of $\widehat{F}$:
$$E^*\frac{\widehat{\omega}}{\widehat{F}}=u\dd y,\quad X=\fraction{g}{u}\fraction{\partial}{\partial y},\quad \coup{\widehat{\mathfrak{I}}_{ \widehat{F}=\infty}}_c=\OO_c .$$
\item $c$ is a singular point of the divisor:
$$E^*\frac{\widehat{\omega}}{\widehat{F}}=\fraction{u}{x^{\epsilon (D_0)}y^{\epsilon (D_1)}}\dd x,\quad X=\fraction{g x^{\epsilon (D_0)}y^{\epsilon (D_1)}}{u}\fraction{\partial}{\partial y},\quad \coup{\widehat{\mathfrak{I}}_{ \widehat{F}=\infty}}_c=\OO_c.$$
\item  $c$ is a zero of $\widehat{F}$:
$$E^*\frac{\widehat{\omega}}{\widehat{F}}=\fraction{u}{y}\dd y,\quad X=\fraction{g}{u}y\fraction{\partial}{\partial y},\quad \coup{\widehat{\mathfrak{I}}_{ \widehat{F}=\infty}}_c=\OO_c.$$
\item $c$ is a pole of $\widehat{F}$:
$$E^*\frac{\widehat{\omega}}{\widehat{F}}=uy\dd y,\quad X=\fraction{g}{uy}\fraction{\partial}{\partial y},\quad \coup{\widehat{\mathfrak{I}}_{ \widehat{F}=\infty}}_c=(y).$$

\end{enumerate}
Hence, in any case, the morphism $E^*\frac{\widehat{\omega}}{\widehat{F}}(\cdot)$ is onto the sheaf $\widehat{\mathfrak{I}}_{ \widehat{F}=\infty}$, which ensures the exactness of the sequence.
\end{demo}

\begin{coro}\label{corobs}
The space of formal infinitesimal obstructions to the problem $(\mathcal{P})$ is $H^1(\DD,\widehat{\mathfrak{I}}_{\widehat{F}=\infty})$. This is a $\mathbb{C}$-space of finite dimension. This dimension is a topological invariant. 
\end{coro}
\noindent In order to prove the above corollary, let us first consider $P$ a germ of curve in $(\mathbb{C}^2,0)$ and $\widehat{\mathfrak{I}}_P$ the subsheaf of $\OO_\MM$ of functions vanishing along the strict transform of $P$ by $E$. The map $E$ is a composition $E=E_0\circ\cdots \circ E_N$ where $E_i$ is the standard blowing-up of one point. For any point $c$ in a divisor $(E_0\circ \cdots\circ E_j)^{-1}(0)$, $0\geq j\geq N$, $\nu_{c}(P)$ refers to the multiplicity at $c$ of the strict transform of $P$ with respect to $E_0\circ\cdots\circ E_j$.
\begin{lemme}
$$\dim_{\mathbb{C}}H^1(\DD,\widehat{\mathfrak{I}}_P)=\sum_{c}\fraction{v_c(P)(v_c(P)-1)}{2}$$
\end{lemme}
\begin{demo} The proof is an induction on the length of the blowing-up process. Let $E_0$ be the blowing-up of the origin.  Let us consider the canonical system of coordinates $(x_1,y_1)$ and $(x_2,y_2)$ in adapted neighborhood of $E_0^{-1}(0)$ such that the change of coordinates is written
$$y_2=y_1x_1, \  x_2=\fraction{1}{y_1}.$$
Let $p$ be a reduced equation of $P$ and $p_1$ and $p_2$ defined by
$$E^*p=x_1^{\nu_0(p)}p_1\quad E^*p=y_2^{\nu_0(p)}p_2.$$
Hence, we can describe the space of global sections

\begin{eqnarray*}
H^0(V_1,\widehat{\mathfrak{I}}_{P}) &\simeq & p_1 \mathbb{C}[[x_1,y_1]] \\
H^0(V_1\cap V_2,\widehat{\mathfrak{I}}_{P}) &\simeq & p_1 \mathbb{C}[[x_1]]((y_1)) \\
H^0(V_2,\widehat{\mathfrak{I}}_{ P}) &\simeq & \left\{\left. y_1^{-\nu_0(p)}\sum_{i,j\in \mathbb{N}^2}a_{ij}x_1^jy_1^{j-i}\right| a_{ij}\in\mathbb{C}\right\} 
\end{eqnarray*}
An simple computation shows that the following isomorphism 
$$\left. { H^0(V_1\cap V_2,\widehat{\mathfrak{I}}_{ P})}\right/{H^0(V_1,\widehat{\mathfrak{I}}_{ P}) \bigoplus H^0(V_2,\widehat{\mathfrak{I}}_{P})}\simeq \mathbb{C}^{\frac{\nu_0(p)(\nu_0(p)-1)}{2}}.$$ 
We decompose the desingularization morphism $E=E_0\circ E_1$ where $E_0$ is the first blowing-up of the origin. Let $\{s_1,\ldots,s_n\}$ refer to the intersection of $D_0$ and the strict transform of $P$. We denote the components of the exceptional divisor $\DD_i=E_1^{-1}(s_i)$. For $i=1,\ldots,n$, consider $U_i(\epsilon)=B(s_i,\epsilon),\ \epsilon >0$ be a disc for any smooth metric on $D_0$ such that $U_i$ does not meet $U_j$ for $j\neq i$. Let $U_0$ be the complementary of $\displaystyle\bigcup_{i=1,\ldots,4}\overline{B}(s_i,\epsilon/2)$. Finally, let us denote by $\UU_i$ the open set $E_1^{-1}(U_i)$.
The system $\{\UU_0,\UU_1,\ldots,\UU_i\}$ provides a covering of the divisor $\DD$ and the Mayer-Vietoris sequence for the sheaf $\widehat{\mathfrak{I}}_{P}$ is written

\begin{equation*}
0\rightarrow N \rightarrow H^1(\DD,\widehat{\mathfrak{I}}_{ P})\rightarrow \bigoplus_i H^1(\UU_i,\widehat{\mathfrak{I}}_{ P})\rightarrow  \bigoplus_{ij} H^1(\UU_i\cap \UU_j,\widehat{\mathfrak{I}}_{ P})\rightarrow 0.
\end{equation*}
where $N$ is given by
$$ \bigoplus_i H^0(\UU_i,\widehat{\mathfrak{I}}_{ P})\rightarrow \bigoplus_{ij} H^0(\UU_i\cap \UU_j,\widehat{\mathfrak{I}}_{ P})\rightarrow N \rightarrow 0.$$
Since $\UU_i\cap \UU_j$ is Stein and 
$\widehat{\mathfrak{I}}_{ P}$ a coherent sheaf, the module $H^1(\UU_i\cap \UU_j,\widehat{\mathfrak{I}}_{ P})$ is trivial. Moreover, the morphism $E_1$ and the Hartogs argument induce following isomorphisms
\begin{eqnarray*}
H^0(\UU_i,\widehat{\mathfrak{I}}_{P})&\simeq & H^0(U_i,\widehat{\mathfrak{I}}_{ P}),\\
H^0(\UU_i\cap \UU_j,\widehat{\mathfrak{I}}_{ P})&\simeq & H^0(U_i\cap U_j,\widehat{\mathfrak{I}}_{ P}).
\end{eqnarray*}
Hence, $N$ is identified with $H^1(D_0,\widehat{\mathfrak{I}}_P)$. All these remarks and an inductive limit on the neighborhood of $\DD_i$ provide the next isomorphisms
$$H^1(\DD,\widehat{\mathfrak{I}}_{ P})\simeq H^1(D_0,\widehat{\mathfrak{I}}_P)\oplus\bigoplus_i H^1(\DD_i,\widehat{\mathfrak{I}}_{ P}).$$
Therefore, the lemma is a straightforward computation from the hypothesis of induction and the formula above.
\end{demo}

\begin{demo}{\it \ (\ref{corobs})} Applying the previous lemma with $P=\{\widehat{F}=\infty\}$, once gets the finite dimension statement. Furthermore, looking at the formula, one can see that this dimension depends only the topology of the morphism of desingularization which is a topological invariant of the foliation \cite{camacho}.
\end{demo}

\subsection{A generic example}
Let $r$ and $n\geq 2$ be positive integers and $\FF_{n,r_1,\ldots,r_n}$ the foliation given by the one form
$$ x^{r+2}\dd \left(\fraction{x^{r+2}-\sum_{j=1}^{r+1}q_jx^{r+1-j}y^j }{x^{r+1}}\right).$$
Let us suppose that $Q(t)=\sum_{j=1}^{r+1}q_jt^j$ satisfies
$Q'(t)=\prod_{i=1}^{n}(t-t_j)^{r_j}$. 
In (\cite{martine}), M. Klughertz proved that these foliations are topological normal forms for $\mathcal{M}$-foliation, i.e., foliations regular after one blowing-up such that all its invariant curves are separatrix. To be more specific, any foliation regular after one blowing-up such that all its invariant curves are separatrix and satisfying
\begin{itemize}
\item there is $n$ invariant curves tangent to the divisor.
\item these $n$ curves are tangent to the divisor with orders $r_1,\ldots, r_n$.
\end{itemize}
is topologically to a foliation $\FF_{n,r_1,\ldots,r_n}$. The position of the tangency point $t_1,\ldots,t_n$ does not have any importance. Since, the dimension of $H^1(\DD,\widehat{\mathfrak{I}}_{\widehat{F}=\infty})$ is a topological invariant, one can compute this dimension for this family of foliations in order to extend the result to any $\MM$-foliation. 

\medskip
\noindent By a direct computation, one has $\nu_0(\FF_{n,r_1,\ldots,r_n})=r+1$. The foliation is regular after one blowing and the exceptionnal divisor $D_0$ is not invariant. There are $n$ integral curves $S_i$, $i=1..n$, which are tangent to the divisor. In the canonical coordinates $y=tx,x=x$, the points of tangency are given by $t=t_j$ and the $r_j$ are the respective order of tangency. Hence, the foliation is completely reduced once one has reduced the germ of curves $S_i\cup(D_0)_{t_i}$. Therefore, the foliation has $n$ isolated separatrix, which are the curves $S_i$ if viewed after one blowing. Let us denote by $h_i$ a reduced equation of $S_i$ at the origin. One can see that $\nu_0(h_i)=r_i+1$. Futhermore, the component $D_0$ is dicritical with valence $n$. Let $\{\alpha_i\}_{i=1..n-2}$ be $n-2$ equations of curves of the pencil of $D_0$. By definition, the meromorphic function
$$F=\fraction{h_1(x,y)h_2(x,y)\cdots h_n(x,y)}{\alpha_1(x,y)\alpha_2(x,y)\cdots \alpha_{n-2}(x,y)}$$
is a balanced equation for $\FF$. In particular, since $\nu_0(\alpha_i)=1$,
$$\nu_0(F)=\sum_{i=1}^n\nu_0(h_i)-(n-2)=\sum_{i=1}^n(r_i+1)-n+2=r+2$$
Hence,
$$\nu_0(F)=\nu_0(\FF_{n,r_1,\ldots,r_n})+1$$
which was predicted by (\ref{critere2esp}). Moreover, in view of (\ref{corobs}), the space of obstructions is of dimension $\fraction{(n-2)(n-3)}{2}$. Hence, for any $\MM$-simple foliation $\FF$ topologically equivalent to one $\FF_{n,r_1,\ldots,r_n}$ with $n$ greater than $4$, there exists a deformation of the resolution space, which cannot carry any isoholonomical deformation of $\FF$. As already mentioned, this situation never occurs in the non-dicritical case and suggests that the space of leaves and the resolution space share some analytical invariants.

{\small 
\bibliography{biblithese}
\bibliographystyle{plain}}

\end{document}